\documentclass[twoside]{article} 
\usepackage{graphicx}
\date{} 
\title{On the $\nu$-zeros of the modified Bessel function $K_{i\nu}(x)$ of positive argument}
\author{\sc R. B.\ Paris \\
{\em Division of Computing and Mathematics,} \\
{\em Abertay University, Dundee DD1 1HG, UK}}
\begin{document}
\def\f#1#2{\mbox{${\textstyle \frac{#1}{#2}}$}}
\def\dfrac#1#2{\displaystyle{\frac{#1}{#2}}}
\def\boldal{\mbox{\boldmath $\alpha$}}
\newcommand{\bee}{\begin{equation}}
\newcommand{\ee}{\end{equation}}
\newcommand{\sa}{\sigma}
\newcommand{\ka}{\kappa}
\newcommand{\al}{\alpha}
\newcommand{\la}{\lambda}
\newcommand{\ga}{\gamma}
\newcommand{\eps}{\epsilon}
\newcommand{\om}{\omega}
\newcommand{\fr}{\frac{1}{2}}
\newcommand{\fs}{\f{1}{2}}
\newcommand{\g}{\Gamma}
\newcommand{\br}{\biggr}
\newcommand{\bl}{\biggl}
\newcommand{\ra}{\rightarrow}
\newcommand{\gtwid}{\raisebox{-.8ex}{\mbox{$\stackrel{\textstyle >}{\sim}$}}}
\newcommand{\ltwid}{\raisebox{-.8ex}{\mbox{$\stackrel{\textstyle <}{\sim}$}}}
\renewcommand{\topfraction}{0.9}
\renewcommand{\bottomfraction}{0.9}
\renewcommand{\textfraction}{0.05}
\newcommand{\mcol}{\multicolumn}
\date{}
\maketitle
\pagestyle{myheadings}
\markboth{\hfill \sc R. B.\ Paris  \hfill}
{\hfill \sc Zeros of the modified Bessel function\hfill}
\begin{abstract}
The modified Bessel function of the second kind $K_{i\nu}(x)$ of imaginary order for fixed $x>0$ possesses a countably infinite sequence of real zeros. Recently it has been shown that the $n$th zero behaves like $\nu_n\sim \pi n/\log\,n$ as $n\to\infty$. In this note we determine a more precise estimate for the bahaviour of these zeros for large $n$ by making use of the known asymptotic expansion of $K_{i\nu}(x)$ for large $\nu$.
Numerical results are presented to illustrate the accuracy of the expansion obtained.
\vspace{0.3cm}

\vspace{0.1cm}
 
\end{abstract}

\vspace{0.3cm}

\noindent $\,$\hrulefill $\,$

\vspace{0.3cm}

\begin{center}
{\bf 1.\ Introduction}
\end{center}
\setcounter{section}{1}
\setcounter{equation}{0}
\renewcommand{\theequation}{\arabic{section}.\arabic{equation}}
The modified Bessel function of the second kind (the Macdonald function) $K_{i\nu}(x)$ of purely imaginary order and argument $x>0$ is given by
\bee\label{e11}
K_{\i\nu}(x)=\int_0^\infty e^{-x\cosh t} \cos \nu t\,dt=\frac{1}{2}\int_{-\infty}^\infty e^{-x\cosh t+i\nu t}dt.
\ee
In a recent paper, Bagirova and Khanmamedov \cite{BK} have shown that $K_{i\nu}(x)$ has a countably infinite number of (simple) real zeros in $\nu$ when $x>0$ is fixed. We label the zeros $\nu_n\equiv\nu_n(x)$ ($n=0, 1, 2 \ldots$) and observe that it is sufficient to consider only the case $\nu>0$ since $K_{i\nu}(x)=K_{-i\nu}(x)$. By transforming the differential equation satisfied by $K_{i\nu}(x)$ into a one-dimensional Schr\"odinger equation with an exponential potential, these authors employed the well-known quantisation rule  to deduce the leading asymptotic behaviour of the $n$th zero given by
\bee\label{e12}
\nu_n\sim \frac{\pi n}{\log\,n}\qquad (n\to+\infty).
\ee

In this note, we consider the behaviour of the large-$n$ $\nu$-zeros of $K_{\i\nu}(x)$ in more detail. To achieve this we make use of the known asymptotic expansion of $K_{i\nu}(x)$ for $\nu\to+\infty$ stated in Section 2. A typical plot of $e^{\pi\nu/2} K_{i\nu}(x)$ as a function of $\nu$ is shown in Fig.~1. 
\begin{figure}[th]
	\begin{center}	\ \includegraphics[width=0.50\textwidth]{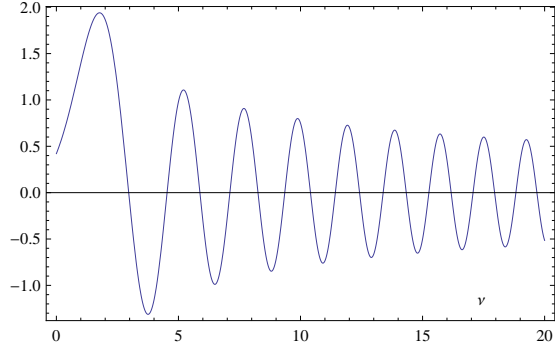}
\caption{\small{Plot of $e^{\pi\nu/2} K_{i\nu}(x)$ as a function of $\nu$ when $x=1$.}}
\end{center}
\end{figure}
\vspace{0.6cm}

\begin{center}
{\bf 2.\ Derivation of the equation describing the large-$n$ zeros}
\end{center}
\setcounter{section}{2}
\setcounter{equation}{0}
\renewcommand{\theequation}{\arabic{section}.\arabic{equation}}
We start with the asymptotic expansion of $K_{i\nu}(x)$ for $\nu\to+\infty$ given by \cite{NMT} (see also\footnote{There is a misprint in (1.4.8) of \cite{PHad}: the argument of the Bessel function should be $x$.}  \cite[pp.~41--42]{PHad})
\bee\label{e21}
e^{\pi\nu/2} K_{i\nu}(x)\sim \Re\,\sqrt{\frac{2\pi}{\nu \tanh \mu}} e^{i\Phi} \sum_{k=0}^\infty \frac{(\fs)_k C_k(\mu)}{(\fs i\nu \tanh \mu)^k},
\ee
where
\[\Phi:=\nu(\mu-\tanh \mu)-\frac{\pi}{4},\quad\cosh \mu=\frac{\nu}{x}\]
and $(a)_k=\g(a+k)/\g(a)$ denotes the Pochhammer symbol. The first few coefficients $C_k(\mu)$ are
\[C_0(\mu)=1,\quad C_1(\mu)=-\frac{1}{8}+\frac{5 \coth^2 \mu}{24},\quad C_2(\mu)=\frac{3}{128}-\frac{77 \coth^2\mu}{576}+\frac{385\coth^4\mu}{3456}~;\]
values of $C_k(\mu)$ for $k\leq 5$ are derived in the appendix. It is important to point out that the determination of the expansion (\ref{e21}) from the integral (\ref{e11}) involves (when $\nu>x$) an infinite number of contributing saddle points. The expansion (\ref{e21}) results from the two dominant saddles, the remaining saddles yielding an exponentially small contribution of O$(e^{-\pi\nu})$.
Then, from (\ref{e21}), the zeros of $K_{i\nu}(x)$ are given asymptotically by
\[\cos \Phi \bl\{1-\frac{3C_2(\mu)}{(\nu\tanh \mu)^2}+\frac{105C_4(\mu)}{(\nu\tanh \mu)^4}+\cdots\br\}\hspace{5cm}\]
\[\hspace{4cm}+\sin \Phi \bl\{\frac{C_1(\mu)}{\nu\tanh \mu}-\frac{15C_3(\mu)}{(\nu\tanh \mu)^3}+\frac{945C_5(\mu)}{(\nu\tanh \mu)^5}+\cdots\br\}=0.\]

Let
\bee\label{e22}
\Phi=(n+\fs)\pi+\epsilon,
\ee
where $n$ is a large positive integer and $\epsilon$ is a small quantity. Then
\[\tan \epsilon={\displaystyle\frac{\frac{C_1(\mu)}{\nu \tanh \mu}-\frac{15C_3(\mu)}{(\nu\tanh \mu)^3}+\frac{945C_5(\mu)}{(\nu\tanh \mu)^5}+\cdots}{1-\frac{3C_2(\mu)}{(\nu\tanh \mu)^2}+\frac{105C^4(\mu)}{(\nu\tanh \mu)^4}+\cdots}}\]
\[=\frac{1}{\nu\tanh \mu}\bl\{C_1(\mu)+\frac{3(C_1(\mu)C_2(\mu)-5C_3(\mu))}{(\nu\tanh \mu)^2}\]
\[+\frac{9C_1(\mu)C_2^2(\mu)-45C_2(\mu)C_3(\mu)-105C_1(\mu)C_4(\mu)+945C_5(\mu)}{(\nu\tanh \mu)^4}+\cdots\br\},
\]
so that using $\arctan z=z-z^3/3+z^5/5- \ldots\,$ ($|z|<1$), we obtain
\[\epsilon=\frac{a_0(\mu)}{\nu}+\frac{a_1(\mu)}{\nu^3}+\frac{a_2(\mu)}{\nu^5}+\cdots,\]
where
\[a_0(\mu)=\frac{C_1(\mu)}{\tanh \mu},\qquad a_1(\mu)=\frac{3(C_1(\mu)C_2(\mu)-5C_3(\mu))-\frac{1}{3}C_1^3(\mu)}{\tanh^3\mu},\]
\[a_2(\mu)=\frac{1}{\tanh^5\mu}\bl\{3(3C_1(\mu)C_2^2(\mu)-15C_2(\mu)C_3(\mu)-35(C_1(\mu)C_4(\mu)-9C_5(\mu)))\]
\[-3C_1^2(\mu)(C_1(\mu)C_2(\mu)-5C_3(\mu))+\frac{1}{5}C_1^5(\mu)\br\}.\]

Since $\tanh \mu=\sqrt{1-x^2/\nu^2}$, the coefficients $a_k(\mu)$ possess expansions in inverse powers of $\nu^2$. Some laborious algebra shows that
\begin{eqnarray*}
a_0(\mu)&=&\frac{1}{12}+\frac{x^2}{4\nu^2}+\frac{11x^4}{32\nu^4}+O(\nu^{-6})\\
a_1(\mu)&=&\frac{1}{360}-\frac{x^2}{4\nu^2}+O(\nu^{-4})\\
a_2(\mu)&=&\frac{1}{1260}+O(\nu^{-2}),
\end{eqnarray*}
whence we obtain
\bee\label{e23}
\epsilon=\frac{1}{12\nu}+\frac{1}{\nu^3}\bl(\frac{1}{360}+\frac{x^2}{4}\br)+\frac{1}{\nu^5}\bl(\frac{1}{1260}-\frac{x^2}{4}+\frac{11x^4}{32}\br)+\cdots.
\ee

Making use of the expansion
\[\mu-\tanh \mu=\log\,\la\nu+\frac{x^2}{4\nu^2}+\frac{x^4}{32\nu^4}+\frac{x^6}{96\nu^6}+ O(\nu^{-8}),\qquad \la:=\frac{2}{ex},\]
we then finally obtain from (\ref{e22}) and (\ref{e23}) the equation describing the large-$n$ zeros of $K_{i\nu}(x)$ given by 
\bee\label{e24}
\nu\log\,\la\nu=m+\frac{A_0}{\nu}+\frac{A_1}{\nu^3}+\frac{A_2}{\nu^5}+\cdots,
\ee
where
\[A_0=\frac{1}{12}-\frac{x^2}{4},\quad A_1=\frac{1}{360}+\frac{x^2}{4}-\frac{x^4}{32},\quad A_2=\frac{1}{1260}-\frac{x^2}{4}+\frac{11x^4}{32}-\frac{x^6}{96}\]
and, for convenience, we have put $m:=(n+\f{3}{4})\pi$.
\vspace{0.6cm}

\begin{center}
{\bf 3.\ Solution of the equation for the zeros}
\end{center}
\setcounter{section}{3}
\setcounter{equation}{0}
\renewcommand{\theequation}{\arabic{section}.\arabic{equation}}
To solve (\ref{e24}) we expand $\nu$ as
\[\nu_n=\xi+\frac{c_0}{\xi}+\frac{c_1}{\xi^3}+\frac{c_2}{\xi^5}+\cdots ,\]
where the $c_k$ are constants to be determined and we suppose that $\xi$ is large as $n\to\infty$. Substitution in (\ref{e24}) then produces
\[\xi \log \la\xi+\frac{c_0(1+\log \la\xi)}{\xi}+\frac{c_1(1+\log \la\xi)+\fs c_0^2}{\xi^3}+\frac{c_2(1+\log \la\xi)+c_0c_1-\f{1}{6}c_0^3}{\xi^5}+\cdots\]
\[=m+\frac{A_0}{\xi}+\frac{A_1-A_0c_0}{\xi^3}+\frac{A_2-3A_1c_0+A_0(c_0^2-c_1)}{\xi^5}+\cdots\,.
\]
Equating coefficients of like powers of $\xi$, we obtain
\bee\label{e31}
\xi\log \la\xi=m,
\ee
and
\[c_0=\frac{A_0}{1+\log \la\xi},\quad c_1=\frac{A_1-A_0 c_0-\fs c_0^2}{1+\log \la\xi},\]
\[c_2=\frac{A_2-3A_1c_0+A_0(c_0^2-c_1)-c_0c_1+\f{1}{6}c_0^3}{1+\log \la\xi}~.\]
The solution of (\ref{e31}) for the lowest-order term $\xi$ can be expressed in terms of the Lambert $W$ function, which is the (positive) solution\footnote{In \cite[p.~111]{DLMF} this is denoted by Wp$(z)$.} of $W(z) e^{W(z)}=z$ for $z>0$.
Rearrangement of (\ref{e31}) shows that
\[\frac{m}{\xi}\,e^{m/\xi}=\la m,\]
whence
\bee\label{e33}
\xi=\frac{m}{W(\la m)}.
\ee

The asymptotic expansion of $W(z)$ for $z\to+\infty$ is \cite{CK}, \cite[(4.13.10)]{DLMF} 
\[W(z)\sim L_1-L_2+\frac{L_2}{L_1}+\frac{L_2^2-2L_2}{2L_1^2}+\frac{(6-9L_2+2L_2^2)L_2}{6L_1^3}+\frac{(-12+36L_2-22L_2^2+3L_2^3)L_2}{12L_1^4}+ \cdots,\]
where $L_1=\log z$, $L_2=\log \log z$, from which it follows that
\[\frac{1}{W(z)}\sim\frac{1}{L_1}\bl\{1+\frac{L_2}{L_1}+\frac{L_2(L_2-1)}{L_1^2}+\frac{(1-5L_2+2L_2^2)L_2}{2L_1^3}+
\frac{(-6+21L_2-26L_2^2+6L_2^3)L_2}{6L_1^4}+\cdots\br\}.\]
Then we have the expansion for $\xi$ as $\la m\to\infty$ given by
\bee\label{e34}
\xi\sim\frac{m}{\log \la m}\bl\{1+\frac{\log\log \la m}{\log \la m}+\frac{\log\log \la m (\log\log \la m-1)}{(\log \la m)^2}+ \cdots\br\},
\ee
where we recall that $m=(n+\f{3}{4})\pi$.

If we define $\chi:=\xi/m$, so that by (\ref{e31}) $1+\log \la\xi=(1+\chi)/\chi$ and introduce the coefficients
\bee\label{e3c}
B_0=\frac{A_0}{1+\chi},\quad B_1=\frac{A_1-A_0 c_0-\fs c_0^2}{\chi^2 (1+\chi)},\quad 
B_2=\frac{A_2-3A_1c_0+A_0(c_0^2-c_1)-c_0c_1+\f{1}{6}c_0^3}{\chi^4 (1+\chi)},
\ee
then we finally have the result:
\newtheorem{theorem}{Theorem}
\begin{theorem}$\!\!\!.$\ \ 
The expansion for the $n$th $\nu$-zero of $K_{\i\nu}(x)$ for fixed $x>0$ is
\bee\label{e32}
\nu_n\sim \frac{m}{W(\la m)}+\frac{B_0}{m}+\frac{B_1}{m^3}+\frac{B_2}{m^5}+\cdots\qquad (n\to\infty),
\ee 
where $m=(n+\f{3}{4})\pi$, $\la=2/(ex)$ and the coefficients $B_k$ are given in (\ref{e3c}).
\end{theorem}
\vspace{0.6cm}

\begin{center}
{\bf 4.\ Numerical results}
\end{center}
\setcounter{section}{4}
\setcounter{equation}{0}
\renewcommand{\theequation}{\arabic{section}.\arabic{equation}}
The leading form of $\xi$ from (\ref{e34}) is, with $\la=2/(ex)$,
\[\xi\sim \frac{(n+\f{3}{4})\pi}{\log (n+\f{3}{4})\pi+\log \la}\qquad(n\to\infty),\]
which yields the approximate estimate in (\ref{e12}). In numerical calculations it is found more expedient to use the expression for $\xi$ in terms of the Lambert function in (\ref{e33}), rather than (\ref{e34}), since the asymptotic scale in this latter series is $\log \la m$ and so requires an extremely large value of $n$ to attain reasonable accuracy.

We present numerical results in Table 1 showing the zeros of $K_{i\nu}(x)$ computed using the FindRoot command in {\it Mathematica} compared with the asymptotic values determined from the expansion (\ref{e32}) with coefficients $B_k$, $k\leq 2$, where $\xi$ is evaluated from (\ref{e33}). The value of the zeroth-order approximation $\xi=m/W(\la m)$ is shown in the final column. It is seen that there is excellent agreement with the computed zeros, even for $n=1$.
\begin{table}[h]
\caption{\footnotesize{Values of the zeros of $K_{i\nu}(x)$ and their asymptotic estimates when $x=1$.}}
\begin{center}
\begin{tabular}{|r|r|r|r|}
\hline
&&&\\[-0.3cm]
\mcol{1}{|c|}{$n$} & \mcol{1}{c|}{$\nu_n$} & \mcol{1}{c|}{Asymptotic} & \mcol{1}{c|}{$\xi$}\\
\hline
&&&\\[-0.3cm]
1  & 4.5344907181 & 4.534{\bf 5}024086 & 4.550063\\
2  & 5.8798671997 & 5.87986{\bf 8}9800 & 5.890918\\
4  & 8.2589364092 & 8.258936{\bf 5}588 & 8.265990\\
5  & 9.3550938258 & 9.3550938{\bf 8}60 & 9.361083\\
10 & 14.3318529171 & 14.33185291{\bf 9}8 & 14.335296\\
15 & 18.8230418511 & 18.823041851{\bf 4} & 18.825473\\
20 & 23.0318794957 & 23.031879495{\bf 8} & 23.033764\\
30 & 30.9169674670 & 30.9169674670 & 30.918273\\
\hline
\end{tabular}
\end{center}
\end{table}

\vspace{0.6cm}

\begin{center}
{\bf Appendix:  The coefficients $C_k(\mu) $ in the expansion (\ref{e21})}
\end{center}
\setcounter{section}{1}
\setcounter{equation}{0}
\renewcommand{\theequation}{\Alph{section}.\arabic{equation}}
The integral in (\ref{e11}) can be cast in the form
\bee\label{a1}
K_{i\nu}(x)=\frac{1}{2}\int_{-\infty}^\infty e^{-x\psi(t)}dt,\qquad \psi(t):=\cosh t-it \cosh \mu,
\ee
where $\cosh \mu=\nu/x$. In the right-half plane, the dominant saddle point (where $\psi'(t)=0$) is situated at $t_0=\mu+\fs\pi i$, with a similar saddle in the left-half plane at $-\mu+\fs\pi i$.
When $\nu>x$, there are in addition two infinite strings of subdominant saddles at $\pm\mu+(2k+\fs)\pi i$ ($k=1, 2, \ldots$) parallel to the positive imaginary axis over which the integration path in (\ref{a1}) is deformed; see \cite{NMT}, \cite[pp.~41--42]{PHad} for details.

We introduce the new variable $w$ by
\[w^2=\psi(t)-\psi(t_0)=\frac{1}{2}i \sinh \mu \,(t-t_0)^2+\frac{1}{6}i \cosh \mu\, (t-t_0)^3+\frac{1}{24}i \sinh \mu\,(t-t_0)^4+\cdots\]
to find with the help of {\it Mathematica} using the InverseSeries command
(essentially Lagrange inversion)
\[t-t_0=e^{-\pi i/4}\sqrt{\frac{2}{\sinh \mu}}\,w+\frac{i \coth \mu}{2\sinh \mu}\,w^2-ie^{-\pi i/4}\sqrt{\frac{2}{\sinh \mu}}\,\frac{(-3+5\coth^2\mu)}{36\sinh \mu}\,w^3+\cdots .\]
Differentiation then yields the expansion
\bee\label{a2}
\frac{dt}{dw}\stackrel{e}{=}e^{-\pi i/4}\sqrt{\frac{2}{\sinh \mu}}\,\sum_{k\geq0} \frac{C_k(\mu)}{(\fs i\sinh \mu)^k}\,w^{2k},
\ee
where $\stackrel{e}{=}$ signifies the inclusion of {\em only the even powers of $w$}, since odd powers will not enter into this calculation. The first few coefficients are found to be:
\[C_0(\mu)=1,\qquad C_1(\mu)=-\frac{1}{8}+\frac{5\coth^2\mu}{24},\]
\[C_2(\mu)=\frac{3}{128}-\frac{77\coth^2\mu}{576}+\frac{385\coth^4\mu}{3456},\] 
\[C_3(\mu)=
-\frac{5}{1024} + \frac{1521 \coth^2\mu}{25600} - \frac{17017 \coth^4\mu}{138240} + \frac{
 17017 \coth^6\mu}{248832},\]
\[C_4(\mu) = 
\frac{35}{32768} - \frac{96833 \coth^2\mu}{4300800} + \frac{
 144001 \coth^4\mu}{1720320} - \frac{1062347 \coth^6\mu}{9953280} + \frac{
 1062347 \coth^8\mu}{23887872},\]
\[C_5(\mu) = 
-\frac{63}{262144} + \frac{67608983 \coth^2\mu}{8670412800} - \frac{
 35840233 \coth^4\mu}{796262400} + \frac{3094663 \coth^6\mu}{31850496}\]
 \bee\label{a3}
  - \frac{
 154040315 \coth^8\mu}{1719926784} + \frac{154040315 \coth^{10}\mu}{5159780352}.
 \ee

Then the integral (\ref{a1}) resulting from the saddle $t_0$ can be expressed as 
\[\int_{-\infty}^\infty e^{-xw^2} \,\frac{dt}{dw}\,dw,\]
which, after substitution of the expansion (\ref{a2}) followed by routine integration and taking into account the contribution from the saddle at $t=-\mu+\fs\pi i$, then yields the expansion for $K_{i\nu}(x)$ stated in (\ref{e21}).

\end{document}